\documentclass[12pt]{amsart}
\usepackage{thmtools}
\usepackage{amsfonts,amsmath,amssymb}
\usepackage{hyperref}
\usepackage{paralist}
\usepackage{stmaryrd}
\usepackage[linesnumbered,ruled,vlined,norelsize]{algorithm2e}
\usepackage[utf8]{inputenc}

\declaretheoremstyle[bodyfont=\normalfont]{noncursive}
\declaretheorem{theorem}
\declaretheorem[numberwithin=section]{lemma}

\declaretheorem[numberlike=lemma]{proposition}
\declaretheorem[numberlike=lemma]{corollary}
\declaretheorem[style=noncursive,numberlike=lemma]{definition}

\declaretheorem[style=noncursive,numberlike=lemma]{remark}

\DeclareMathOperator{\ord}{ord}

\usepackage{amsmath,amssymb,amsthm,enumerate}

\begin{document}
\numberwithin{equation}{section}

\def\1#1{\overline{#1}}
\def\2#1{\widetilde{#1}}
\def\3#1{\widehat{#1}}
\def\4#1{\mathbb{#1}}
\def\5#1{\frak{#1}}
\def\6#1{{\mathcal{#1}}}

\def\C{{\4C}}
\def\R{{\4R}}
\def\N{{\4N}}
\def\Z{{\4Z}}

\def \im{\text{\rm Im }}
\def \re{\text{\rm Re }}
\def \Char{\text{\rm Char }}
\def \supp{\text{\rm supp }}
\def \codim{\text{\rm codim }}
\def \Ht{\text{\rm ht }}
\def \Dt{\text{\rm dt }}
\def \hO{\widehat{\mathcal O}}
\def \cl{\text{\rm cl }}
\def \bR{\mathbb R}
\def \bC{\mathbb C}
\def \bP{\mathbb P}
\def \C{\mathbb C}
\def \bL{\mathbb L}
\def \bZ{\mathbb Z}
\def \bN{\mathbb N}
\def \scrF{\mathcal F}
\def \scrK{\mathcal K}
\def \scrM{\mathcal M}
\def \cR{\mathcal R}
\def \scrJ{\mathcal J}
\def \scrA{\mathcal A}
\def \scrO{\mathcal O}
\def \scrV{\mathcal V}
\def \scrL{\mathcal L}
\def \scrE{\mathcal E}

\oddsidemargin=0.1in \evensidemargin=0.1in \textwidth=6.4in
\headheight=.2in \headsep=0.1in \textheight=8.4in
\newcommand{\rl}{{\mathbb{R}}}
\newcommand{\cx}{{\mathbb{C}}}
\newcommand{\id}{{\mathbb{I}}}
\newcommand{\m}{{\mathcal{M}}}
\newcommand{\dbar}{\overline{\partial}}
\newcommand{\Db}[1]{\frac{\partial{#1}}{\partial\overline{z}}}
\newcommand{\abs}[1]{\left|{#1}\right|}
\newcommand{\e}{\varepsilon}
\newcommand{\tmop}[1]{\ensuremath{\operatorname{#1}}}
\renewcommand{\Re}{\tmop{Re}}
\renewcommand{\Im}{\tmop{Im}}
\newcommand{\dist}{{\mathrm{dist}}}
\newcommand {\OO}{{\mathcal O}}
\newcommand {\Sig}{{\Sigma}}
\renewcommand {\a}{\alpha}
\newcommand {\MC}{M^{\mathbb C}}
\renewcommand {\b}{\beta}
\newcommand {\Q}{\mathcal Q}
\newcommand{\dop}[1]{\frac{\partial}{\partial #1}}
\newcommand{\dopt}[2]{\frac{\partial #1}{\partial #2}}
\newcommand{\vardop}[3]{\frac{\partial^{|#3|} #1}{\partial {#2}^{#3}}}

\newcommand{\CC}[1]{\mathbb{C}^{#1}}
\newcommand{\CP}[1]{\mathbb{CP}^{#1}}
\newcommand{\RR}[1]{\mathbb{R}^{#1}}
\newcommand{\dw}{\frac{\partial}{\partial w}}
\newcommand{\dz}{\frac{\partial}{\partial z}}
\numberwithin{equation}{section}
\newcommand{\diffcr}[1]{\rm{Diff}_{CR}^{#1}}
\newcommand{\Hol}[1]{\rm{Hol}^{#1}}
\newcommand{\Aut}[1]{\rm{Aut}^{#1}}
\newcommand{\hol}[1]{\mathfrak{hol}^{#1}}
\newcommand{\aut}[1]{\mathfrak{aut}^{#1}}
\newcommand{\fps}[1]{\C\llbracket #1 \rrbracket}
\newcommand{\fpstwo}[2]{#1\llbracket #2 \rrbracket}
\newcommand{\cps}[1]{\C\{#1\}}

\title[]{Regularity of CR-mappings between Fuchsian type hypersurfaces in $\CC{2}$}

\author{Peter Ebenfelt}
\address{Department of Mathematics, University of California at San Diego, La Jolla, CA 92093-0112}
\email{pebenfelt@ucsd.edu}

\author{Ilya Kossovskiy}
\address{Department of Mathematics, Masaryk University, Brno, Czechia//
Department of Mathematics, University of Vienna, Vienna, Austria}
\email{kossovskiyi@math.muni.cz}

\author{Bernhard Lamel}
\address{Department of Mathematics, University of Vienna, Vienna, Austria}
\email{bernhard.lamel@univie.ac.at}


\thanks{The first author was supported in part by the NSF grant DMS-1600701. The second author was supported in part by the Czech Grant Agency (GACR) and the Austrian Science Fund (FWF). The third author
was supported in part by the Austrian Science Fund (FWF)}



\begin{abstract}
We investigate regularity of CR-mappings between real-analytic infinite type hypersurfaces in $\CC{2}$.  We show that, under the {\em Fuchsian type} condition, all (respectively formal or smooth) CR-diffeomorphisms between them are automatically analytic. The Fuchsian condition appears to be in a certain sense optimal for the regularity problem. 
\end{abstract}

\maketitle   

\def\Label#1{}


\def\cn{{\C^n}}
\def\cnn{{\C^{n'}}}
\def\ocn{\2{\C^n}}
\def\ocnn{\2{\C^{n'}}}


\def\dist{{\rm dist}}
\def\const{{\rm const}}
\def\rk{{\rm rank\,}}
\def\id{{\sf id}}
\def\aut{{\sf aut}}
\def\Aut{{\sf Aut}}
\def\CR{{\rm CR}}
\def\GL{{\sf GL}}
\def\Re{{\sf Re}\,}
\def\Im{{\sf Im}\,}
\def\span{\text{\rm span}}

\def\codim{{\rm codim}}
\def\crd{\dim_{{\rm CR}}}
\def\crc{{\rm codim_{CR}}}

\def\lr{\longrightarrow}
\def\phi{\varphi}
\def\eps{\varepsilon}
\def\d{\partial}
\def\a{\alpha}
\def\b{\beta}
\def\g{\gamma}
\def\G{\Gamma}
\def\D{\Delta}
\def\Om{\Omega}
\def\k{\kappa}
\def\l{\lambda}
\def\L{\Lambda}
\def\z{{\bar z}}
\def\w{{\bar w}}
\def\Z{{\1Z}}
\def\t{\tau}
\def\th{\theta}

\def\H{\hat H}

\emergencystretch15pt
\frenchspacing

\newtheorem{Thm}{Theorem}[section]
\newtheorem{Cor}[Thm]{Corollary}
\newtheorem{Pro}[Thm]{Proposition}
\newtheorem{Lem}[Thm]{Lemma}

\theoremstyle{definition}\newtheorem{Def}[Thm]{Definition}

\theoremstyle{remark}
\newtheorem{Rem}[Thm]{Remark}
\newtheorem{Exa}[Thm]{Example}
\newtheorem{Exs}[Thm]{Examples}

\def\bl{\begin{Lem}}
\def\bl{\begin{Lem}}
\def\el{\end{Lem}}
\def\bp{\begin{Pro}}
\def\ep{\end{Pro}}
\def\bt{\begin{Thm}}
\def\et{\end{Thm}}
\def\bc{\begin{Cor}}
\def\ec{\end{Cor}}
\def\bd{\begin{Def}}
\def\ed{\end{Def}}
\def\br{\begin{Rem}}
\def\er{\end{Rem}}
\def\be{\begin{Exa}}
\def\ee{\end{Exa}}
\def\bpf{\begin{proof}}
\def\epf{\end{proof}}
\def\ben{\begin{enumerate}}
\def\een{\end{enumerate}}
\def\beq{\begin{equation}}
\def\eeq{\end{equation}}

\bigskip

\begin{center}\em
Dedicated to the memory of Nick Hanges
\end{center}

\bigskip

\mbox{}

\bigskip


\section{Introduction}

The problem of regularity of CR-maps between CR-submanifolds in complex space is of fundamental importance in the field of Several Complex Variables. Starting from the classical work of Cartan \cite{cartan},  Chern and Moser \cite{CM74}, Pinchuk \cite{pinchuksib}, and Lewy \cite{lewy}, a large amount of publications is dedicated to various positive results in this well developed  direction. In particular, when both the source and the target are real-analytic, the expected regularity of smooth CR-maps is $C^{\omega}$, i.e., they are {\em analytic} (this property implies that the CR-maps extend holomorphically to a neighborhood of the source manifold).  We refer the reader to the book of  Baouendi-Ebenfelt-Rothschild \cite{ber}, the survey of Forstneri\'c \cite{forstneric}, the book of Berhano, Cordaro and Hounie \cite{cordaro}, and the introduction in \cite{nonanalytic} for the set-up of the theory of CR-maps, a historic outline of the analyticity problem, its connections with the boundary regularity of holomorphic maps\,/\, the reflection principle, and the connections of the problem to  the theory of linear PDEs.   

In the particularly well studied case of real-analytic hypersurfaces in $\CC{2}$, it has been known for some time that CR-diffeomorphisms of finite D'Angelo type hypersurface are automatically analytic (see e.g. Baouendi-Jacobowitz-Treves \cite{bjt}). (Note that in the $\CC{2}$-case finite D'Angelo type is equivalent to the Hörmander-Kohn bracket-generating condition and Tumanov nonminimality). In the case of infinite type but Levi-nonflat hypersurfaces, when there exists a complex variety $X\subset M$ passing through the reference point $p$ in the source hypersurface $M$, some partial analyticity results are available. For instance, analyticity has been established by Ebenfelt \cite{1nonminimal} for so-called {\em $1$-nonminimal hypersurfaces} (see the notion of {\em nonminimality order} below), and by Ebenfelt-Huang \cite{eh} for the case of maps admitting a one-sided holomorphic extension.           

On the other hand,  in the recent paper \cite{nonanalytic}, Kossovskiy and Lamel discovered the existence of real-analytic hypersurfaces in $\CC{N},\,N\geq 2$ which are $C^\infty$ CR-equivalent, but are inequivalent analytically. In particular, it follows that $C^\infty$ CR-diffeomorphisms between real-analytic Levi-nonflat hypersurfaces in $\CC{2}$ are {\em not} analytic in general. Moreover, it shows that the equivalence problem for nonminimal real-analytic CR-structures is of a more {\em intrinsic} nature, as a map realizing an equivalence does not necessarily arise from the biholomorphic equivalence of the CR-manifolds as submanifolds in complex space. 

A natural question  immediately raised by the results in \cite{nonanalytic} is to identify an optimal class of ``regular'' real-analytic hypersurfaces, for which CR-diffeomorphism are still analytic. 
The goal of the current paper is to address this question in the $\CC{2}$-case. We consider the class of {\em Fuchsian type hypersurfaces} introduced by the authors in \cite{ekl} (this condition is described explicitly in terms of the defining function of a hypersurface), and prove that CR-diffeomorphisms of Fuchsian type hypersurfaces are automatically analytic. We also show the invariance and optimality of the Fuchsian type condition. 

Another result of us concerns the problem of convergence of {\em formal} CR-maps. Similarly to the analyticity issue, this problem has attracted a lot of attention of experts in complex analysis in the last few decades (see, e.g., the survey \cite{lmsurvey} of Lamel-Mir). \autoref{theor5} below establishes a convergence result for formal CR-maps in the Fuchsian type case.

\smallskip

We now formulate the results below in detail. We start with describing the precise class of hypersurfaces considered in this paper. 
In light of the  above, we deal with
 germs of Levi-nonflat real-analytic hypersurfaces $M\subset\CC{2}$
considered near a point of {\em infinite} type $p\in
M$. If $M$ is such a hypersurface, there is a unique germ of a
 complex hypersurface (complex curve)
 $X\subset M$ passing through $p$. The complex hypersurface
 $X$ consists
of all infinite type points in $M$ near $p$, it is nonsingular and we will
 also refer to it as the {\em infinite type locus of $M$}.
We say that $(M,p)$ is of {\em generic infinite type} if the canonical
extension of the  Levi
form
$$\mathcal L_p:\,\,T^{1,0}_p\times T^{1,0}_p\lr \CC{}T_pM/\CC{}T^{\CC{}}_pM$$
 from $M$ to its {\em complexification}
$M^{\CC{}}\subset\CC{2}\times\overline{\CC{2}}$ \,\,locally vanishes
only on the  complexification
$X^{\CC{}}\subset\CC{2}\times\overline{\CC{2}}$ of $X$. (We
refer the reader to Section 2 for
details). If $M$ is a Levi-nonflat real-analytic hypersurface
with infinite type locus $X$, then $M$ must be of
generic infinite type at points $p$ lying outside of a
 proper real-analytic subset of $X$. 

We say that  local holomorphic coordinates $(z,w)$, where $w = u+iv$,   near $p$
are {\em admissible} (for $M$)   if in these coordinates,   $p$ becomes the
origin and   $M$ is given by  \begin{equation}\Label{madmissiblereal}\label{madmissiblereal}
v=\frac{1}{2}u^m\left(\epsilon |z|^2+\sum_{k,l\geq 2}h_{kl}(u)z^k\bar
z^l\right)=: h(z,\bar z,u) ,\quad \epsilon=\pm 1 \end{equation} (such admissible
coordinates always exist under the generic infinite type assumption, see \cite{nonminimalODE}); in
particular, in these coordinates $X = \{ w = 0\}$.  The integer $m\geq 1$ is an
important invariant of an infinite type hypersurface called {\em the nonminimality order}, and $M$ with such nonminimality order is called {\em $m$-nonminimal}. For an even $m$, we can
further normalize $\epsilon$ to be equal to $1$, while for an odd $m$,
$\epsilon$ is a biholomorphic invariant. Note that the form
\eqref{madmissiblereal} is stable under the group of dilations
\begin{equation}\Label{dilations}\label{dilations}  z\mapsto \lambda z, \quad w\mapsto \mu w,
\quad \mu^{1-m}=\epsilon|\lambda|^2, \quad \lambda\in\CC{}\setminus\{0\}, \quad
\mu\in\RR{}. 
\end{equation}


\smallskip

We are now able to describe the Fuchsian condition.
\begin{definition}\Label{fuchsianM}\label{fuchsianM}
An infinite type hypersurface \eqref{madmissiblereal} is called {\em a hypersurface of Fuchsian type}, if its  defining function $h(z,\bar z,u)$ satisfies
\begin{equation}\Label{FM}\label{FM}
\begin{aligned}
&\ord\,h_{22}(w)\geq m-1; \,\ord\,h_{23}(w)\geq 2m-2; \,\ord\,h_{33}(w)\geq 2m-2;\\
&\ord\,h_{2l}(w)\geq 2m-l+2, \,\,4\leq l\leq 2m+1;\\
&\ord\,h_{kl}(w)\geq 2m-k-l+5,\,\,k\geq 3,\,l\geq 3,\,7\leq k+l\leq 2m+4.
\end{aligned}
\end{equation}
\end{definition}
 We point out that

 \smallskip

 \begin{itemize}

\item \noindent The Fuchsian condition requires vanishing of an appropriate part of the $(2m+4)$-jet of the defining function $h$ at $0$;

 \smallskip

\item \noindent  It is easy to see from \ref{FM} that {for $m=1$ the Fuchsian type condition holds automatically, while for $m>1$ it fails to hold in general};

 \smallskip

 \item  \noindent  As will be shown in Section 3, the Fuchsian type property is holomorphically invariant.

   \end{itemize}

   \smallskip

\begin{remark}\Label{comparefuchs}\label{comparefuchs}
The property of being Fuchsian extends earlier versions of this property given  respectively in the work  \cite{nonminimalODE} Kossovskiy-Shafikov, and the work  \cite{analytic} of Kossovskiy-Lamel. In the paper \cite{nonminimalODE}, a Fuchsian property of {\em generically spherical} hypersurfaces \eqref{madmissiblereal} was introduced. It is possible to check that for a generically spherical hypersurface the two notions of being Fuchsian coincide. In the paper \cite{analytic}, general hypersurfaces \eqref{madmissiblereal} were considered, but the notion of Fuchsian type considered there is weaker than that given in \cite{ ekl} and in the present paper; it serves to guarantee the regularity of infinitesimal CR-automorphisms, while the property \eqref{FM} guarantees regularity of {\em arbitrary} CR-maps. The property introduced in \cite{analytic} is more appropriately addressed as {\em weak Fuchsian type}, while the property \eqref{FM} as the (actual) Fuchsian type.
\end{remark}

Now our main analyticity results are as follows.

\begin{theorem}\Label{theor5}\label{theor4}
Let $M,M^*\subset\CC{2}$ be real-analytic  hypersurfaces, and let $M$ be   of Fuchsian type at a point $p\in M$. Let $U$ be an open neighborhood of $p$ in $\CC{2}$. Then any $C^\infty$ CR-diffeomorphism $H:\,M\cap U\longrightarrow M^*$ is analytic.
\end{theorem}

By applying the Hanges-Treves propagation principle \cite{hanges}, we are able to address the regularity at an arbitrary infinite type point.

\begin{theorem}\Label{theor4}\label{theor4'}
Let $M,M^*\subset\CC{2}$ be real-analytic Levi-nonflat hypersurfaces, and  $U$  an open neighborhood of $p$ in $\CC{2}$. Assume that $U\cap M$ contains an Fuchsian type point $q$. Then any $C^\infty$ CR-diffeomorphism $H:\,M\cap U\longrightarrow M^*$ is analytic.

\end{theorem}

We further obtain a result on the convergence of formal power series maps between Fuchsian type hypersurfaces.

\begin{theorem}\Label{theor4}\label{theor5}
Let $M,M^*\subset\CC{2}$ be real-analytic  hypersurfaces, and let $M$ be   of Fuchsian type at a point $p\in M$. Then any formal invertible power series map $H:\,(M,p)\longrightarrow (M^*,p^*),\,p^*\in M^*$ is convergent.
\end{theorem}

\begin{remark}
As follows from the invariance of the Fuchsian type property under formal power series transformations  (see \autoref{fuchsinv} below), the target hypersurface $M^*$ is also of Fuchsian type at the respective point $p^*=H(p)$. 
\end{remark}

\autoref{theor5} extends earlier results in this direction obtained in \cite{JLmz} in the case $m=1$. It also extends, in a certain sense, the result in \cite{analytic} on the regularity of infinitesimal CR-automorphisms of Fuchsian type hypersurfaces to the case of general maps (not necessarily appearing as flows of infinitesimal CR-automorphisms). However, as discussed above, the Fuchsian type condition in \cite{analytic} is more mild and involves only vanishing conditions on the coefficient functions $h_{kl},\,k+l\leq 7$ (unlike the conditions in \eqref{FM}). As arguments in Section 3 below show, the case of a general CR-mapping requires considering {\em all} the coefficients  $h_{kl}$ in \eqref{FM}, as they appear in the complete (singular) system of ODEs determining a CR-map.




\section{Preliminaries}

\subsection{Infinite type real hypersurfaces}\label{ss:nonminimal} 
We recall that if  $M\subset \CC{2}$ is a real-analytic  hypersurface, then
for any $p\in M$ there exist so-called {\em normal coordinates} $(z,w)$ centered
at $p$ for $M$. The coordinates
being normal means that $(z,w)$ is
a local holomorphic coordinate system near $p$ in
which $p=0$ and for which near $0$,
$M$ is defined by an equation of the form
\[ v = F (z, \bar z , u)\]
for some germ $F$ of a holomorphic function on $\CC{3}$ which satisfies
the normality condition
\[ F(z,0,u) = F(0,\bar z, u) = 0\]
and the
{ reality condition} $F(z,\bar z,u) \in \RR{} $ for $(z,u)\in \CC{} \times\RR{}$ close to $0$ (see e.g. \cite{ber}). Equivalently, $v = F(z,\bar z,u)$ defines a real hypersurface,
and in the coordinates $(z,w)$, we have
$Q_{(0,u)} = \left\{ (0,w) \in U\colon w = u \right\}$

We also recall that $M$ is {\em of infinite type} at $p$ if there exists a germ of a nontrivial
complex curve $X\subset M$ through $p$. It turns out that in normal coordinates, such
a curve $X$ is necessarily defined by $w = 0$ (because $X = Q_{0} = \{ w= 0\}$); 
in particular, any such $X$ is nonsingular. 
It also turns out that $M$ is Levi-flat if and only if in normal
coordinates, it is defined by $v = 0$.
Thus a Levi-nonflat
real-analytic hypersurface $M$ is of infinite type at $p$ if and only if in
 normal coordinates $(z,w)$ as above, the
defining function $F$ satisfies  $F(z,\bar z,0) = 0$.
In other words, $M$ is of infinite type if and only if it
 can defined by an equation of the form
\begin{equation}\label{mnonminimal}
v=u^m\psi(z,\bar z,u),
\text{ with } \psi(z,0,u) = \psi(0,\bar z, u)= 0 \text{ and } \psi(z,\bar z,0)\not\equiv 0,
\end{equation}
where $m\geq 1$. 
It turns out that the integer $m\geq 1$ is independent both
the choice of $p\in X$ and also of the choice of
normal coordinates for $M$ at $p$  (see \cite{meylan}), and we
 say that
$M$ is {\em $m$-infinite type} along $X$ (or at $p$).

We are going to utilize a number of different ways to write a defining function.
Throughout this paper, we use the {\em complex defining function} $\Theta$ in
which $M$ is defined by  \[  w = \Theta (z,\bar z,\bar w) ;\] it is obtained
from $F$ by solving the equation  \[ \frac{w - \bar w}{2i} = F \left(z,\bar z,
\frac{w+\bar w}{2} \right) \] for $w$, and it agrees with the function 
defining the  Segre varieties in those coordinates, that is,  $Q_Z = \{ (z,
\Theta(z,\bar Z)) \colon z\in U^z \}$. We are going to make
extensive use of the Segre varieties and refer the reader to
\cite{ber} for a discussion
of their properties in the general case,
and to \cite{nonanalytic} for
specific properties in the infinite type setting.

The complex defining function (in
normal coordinates) satisfies the
conditions  \[ \Theta(z,0,\tau) = \Theta (0, \chi, \tau ) = \tau, \quad
\Theta(z,\chi,\bar \Theta (\chi, z, w)) = w. \] If $M$ is of
$m$-infinite type at $p$, then
$\Theta (z,\chi,\tau) = \tau \theta (z,\chi,\tau)$ and  thus $M$ is defined by
the equation $ w = \bar w \theta(z,\bar z,\bar w) = \bar w  +\bar w^{m}\tilde
\theta(z,\bar z,\bar w)$, where $\tilde \theta $ satisfies
$ \tilde \theta (z,0,\tau) = \tilde
\theta (0,\chi, \tau) = 0$  and   $ \tilde\theta (z,\chi,0)\neq 0$.

We also note  that the {external complexification} $M^{\CC{}}$ of $M$, which
is the hypersurface in $\CC{2}\times\overline{\CC{2}}$ defined by
$M^{\CC{}} = \left\{ (Z, \zeta) \in U\times \bar U \colon Z \in Q_{\bar \zeta} \right\} $,
is conveniently defined as the graph of the complex defining function $\Theta$,
i.e.
$$w= \Theta(z,\chi,\tau).$$


We also introduce the real line
\begin{equation}\label{Gamma}
\Gamma = \{(z,w)\in M \colon z = 0 \} = \{(0,u) \in M \colon u\in\RR{}\}\subset M,
\end{equation}
and recall that
\[ Q_{(0,u)} = \{w = u\}, \quad (0,u)\in \Gamma \]
for $u\in\RR{}$. This  property, as already mentioned, is
actually equivalent to the normality of
the coordinates $(z,w)$. More precisely, for any real-analytic curve
$\gamma$ through $p$ one can find normal coordinates $(z,w)$ in
which a small piece of $\gamma$ corresponds to
 $\Gamma$ in \eqref{Gamma}.

We finally notice that  a real-analytic Levi-nonflat hypersurface
$M\subset\CC{2}$ has infinite type points of two kinds,
which we will
refer to as  {\em generic} and {\em exceptional} infinite type points,
respectively. A generic point $p\in M$ is characterized by the condition that
the complexified Levi form of $M$ only degenerates
on the complexified infinite type locus $w=\tau=0$ near $p$. (The
complexified Levi form is defined similarly to the classical Levi form, but
instead the $(1,0)$ and the $(0,1)$ vector fields are considered on the
complexification $M^{\CC{}}$, see e.g. \cite{ber}). We refer to a
non-generic point  $p$ as {\em exceptional}. We note that the
set of exceptional points is a proper real-analytic subvariety of $X$ and
that $p\in X$ is generic if and only if the Levi-determinant of $M$ vanishes
to order $m$ along {\em any} real curve $\gamma$ passing through $p$ which
is transverse to $X$ at $p$.

A generic infinite type point is characterized in normal coordinates by
requiring in addition to \eqref{mnonminimal}  the
condition $\psi_{z\bar z}(0,0,0,)\neq 0.$ If $p$ is a generic
infinite type point, we can further
simplify $M$ to the form \eqref{madmissiblereal} above, or alternatively to the
{\em exponential form}
\begin{equation}\Label{exponential}\label{exponential}  w=\bar w e^{i\bar
w^{m-1}\varphi(z,\bar z,\bar w)},\quad\mbox{where}\quad\varphi(z,\bar z,\bar
w)=\pm z\bar z+\sum_{k,l\geq 2}\varphi_{kl}(\bar w)z^k{\bar z}^l \end{equation}
(see, e.g., \cite{nonminimalODE}).

\subsection{Real hypersurfaces and second order differential equations.}\label{sub:realhyp2ndorderequ}
There is a natural way to associate to
a  Levi nondegenerate real hypersurface
$M\subset\CC{N}$ a system of second order
holomorphic PDEs with $1$ dependent and $N-1$ independent
variables by using
the Segre family of the hypersurface $M$.
This remarkable construction
 goes back to
E.~Cartan \cite{cartan} and Segre \cite{segre}
 (see also a remark by Webster \cite{webster}),
and was recently revisited in the work of Sukhov
\cite{sukhov1,sukhov2} in the nondegenerate setting, and in the work of Kossovskiy, Lamel and Shafikov  in the degenerate setting
(see\cite{divergence,nonminimalODE,nonanalytic,analytic}). 
For the convenience of the reader, we recall this procedure in the case
$N=2$, but refer to the above references for more details.

So assume that $M\subset\CC{2}$ is
a smooth real-analytic hypersurface passing through the origin
 and $U = U^z \times U^w$ is chosen small enough.
 The  second order holomorphic ODE associated to $M$
 is uniquely determined by the condition that for every
 $\zeta \in U$,  the
 function $h(z,\zeta) = w(z)$ defining the Segre variety
 $Q_\zeta$ as a graph is a solution of this ODE. 
To be more precise, one can show that the Levi-nondegeneracy
of  $M$ (at $0$) implies that
near the origin, the Segre map
$\zeta\mapsto Q_\zeta$ is injective and the Segre family has
the so-called transversality property: if two distinct Segre
varieties intersect at a point $q\in U$, then their intersection
at $q$ is transverse (actually it turns out
that, again due to the Levi-nondegeneracy of $M$,
the Segre varieties passing through a point $p$
are uniquely determined by their tangent spaces $T_p Q_\zeta$).
Thus, $\{Q_\zeta\}_{\zeta\in U}$ is a
2-parameter  family of holomorphic
curves in $U$ with the transversality property, depending
holomorphically on $\bar\zeta$. It follows from
the holomorphic version of the fundamental ODE theorem (see, e.g.,
\cite{MR2363178}) that there exists a unique second order
holomorphic ODE $w''=\Phi(z,w,w')$ such that for each $\zeta \in U$,
$w(z) = h(z, \bar \zeta)$ is one of its solutions.

We can carry out the construction of this ODE concretely by
utilizing the complex defining equation $w=\Theta(z,\chi,\tau)$
introduced above.
Recall that  the Segre
variety $Q_\zeta$ of a point $\zeta=(a,b)\in U$ is  now given
as the graph
\begin{equation} \label{segre0}w (z)=\rho(z,\bar a,\bar b). \end{equation}
Differentiating \eqref{segre0} once, we obtain
\begin{equation}\label{segreder} w'=\rho_z(z,\bar a,\bar b). \end{equation}
The system of equations  \eqref{segre0} and \eqref{segreder}
can be solved, using the implicit function theorem, for
 $\bar a$ and $\bar b$. This gives us holomorphic functions
 $A $ and $ B$ such that
$$
\bar a=A(z,w,w'),\,\bar b=B(z,w,w').
$$
The application of the implicit function theorem is possible
since the
Jacobian of the system consisting of  \eqref{segre0} and \eqref{segreder} with
respect to $\bar a$ and $\bar b$ is just the Levi determinant of $M$
for $(z,w)\in M$ (\cite{ber}). Differentiating \eqref{segreder} once more,
we can substitute  $\bar a = A(z,w,w') $ and $ \bar b=B(z,w,w')$ to obtain
\begin{equation}\label{segreder2}
w''=\rho_{zz}(z,A(z,w,w'),B(z,w,w'))=:\Phi(z,w,w').
\end{equation}
Now \eqref{segreder2} is a holomorphic second order ODE, for which all
of the functions $w(z) = h(z, \zeta)$ are solutions by construction.
We will denote this associated second order ODE by
$\mathcal E = \mathcal{E}(M) $.

More generally it is possible to  associate  a completely integrable PDE
to any of a wide range of
CR-submanifolds (see \cite{sukhov1,sukhov2}) such that the
correspondence $M\to \mathcal E(M)$ has the following fundamental
properties:

\begin{compactenum}[(1)]

\item Every local holomorphic equivalence $F:\, (M,0)\to (M',0)$
between CR-submanifolds is an equivalence between the
corresponding PDE systems $\mathcal E(M),\mathcal E(M')$;

\item The complexification of the infinitesimal automorphism algebra
$\mathfrak{hol}^\omega(M,0)$ of $M$ at the origin coincides with the Lie
symmetry algebra  of the associated PDE system $\mathcal E(M)$
(see, e.g., \cite{olver} for the details of the concept).
\end{compactenum}

In contrast to the case of a finite type real hypersurface described above,
if
$M\subset\CC{2}$ is of infinite type at the origin one cannot associate
to $M$ a regular second order ODE or even a more general PDE system near
the origin such that the Segre varieties are graphs of solutions.
However, in \cite{nonminimalODE} and \cite{analytic},
Kossovskiy, Lamel and Shafikov
found an injective correspondence associating to a  hypersurface
$M\subset\CC{2}$ at a  generic infinite type point a
certain {\em singular} complex ODE $\mathcal E(M)$ with an
isolated  singularity at the origin.
We are going to base our normal form construction on
this construction, which  is therefore extensively used in the paper
(more details are given in Section 3).

We finally point out that at {\em exceptional} infinite type points, one
can still associate
 a system of singular complex ODEs
to a real-analytic hypersurface $M\subset\CC{2}$ (although possibly of
 higher order $k\geq 2$) as in the paper \cite{KLS}
 Kossovskiy-Lamel-Stolovitch.

\subsection{Complex differential equations with an isolated singularity} \label{ss:isolated}\mbox{}
We will again just gather the facts from the classical theory of singular (complex)
differential equations, and
refer  the reader to e.g. \cite{MR2363178}, \cite{vazow},\cite{laine} for
any details.

A linear system $\mathcal L$ of (holomorphic) first order ODEs  on a domain
$G\subset\CC{}$ (or simply a {\em linear system} in a domain $G$)
is an equation of the form
$y'(x)=A(x)y(x)$, where $A \colon G \to \C^{n\times n}$
is a matrix-valued
holomorphic map on $G$
 and $y(x)=(y_1(x),...,y_n(x))$
is an $n$-tuple of (unknown) functions. The set of  solutions of $\mathcal L$
near a point $p\in G$ is isomorphic to $\CC{n}$ by $y \mapsto y(p)$.
Because every germ $y$  of a solution of $\mathcal L$ at $p\in G$
extends analytically along any path $\gamma\subset G$ starting
at~$p$,  any solution $y(x)$ of $\mathcal L$ is defined
 in all of $G$ as a (possibly multi-valued) analytic function.
If $G$ is a punctured disc, centered at $0$, we say
that $\mathcal L$  has an {\em isolated singularity} (at $x=0$).
If $A(x)$ has a pole at the isolated singularity $x=0$,
we say that the system has a {\em meromorphic singularity}. As
the solutions of $\mathcal L$ are holomorphic in any proper sector
$S\subset G$ of a sufficiently small radius with vertex at
$x=0$, it is important to study the behaviour of the solutions as
$x\rightarrow 0$. If for every sector
$S = \{ x \in G \colon |x|<\delta, \alpha < \arg x < \beta  \} $
there exist constants $C>0$ and  $a\in \mathbb R$ such that
for every solution
$y$ of $\mathcal L$ defined in $S$ we have that
$||y(x)||\leq C|x|^a$  holds for $x\in S$, then we say that
$x=0$ is a {\em regular singularity}, otherwise we say it is
 an {\em irregular singularity}.

An important condition ensuring regularity of a singularity is due to
L.~Fuchs: We say that the singular point
$x=0$ is {\em Fuchsian} if $A(x)$ has a pole of order at most
 $ 1$ at $x = 0$. If $0$ is a Fuchsian singularity, then
 $x = 0$ is a  regular
singular point. Another important property of Fuchsian
singularities is that every formal power series solution (at $x=0$)
of the equation is actually  {\em convergent}.
The dynamical system associated
to a Fuchsian singularity corresponds to the dynamical system of the
 vector field $$x\frac{\partial}{\partial x}+A(x)y\frac{\partial}{\partial y},$$
 which is ``almost''\, non-resonant in the sense of Poincar\'e-Dulac.

However, in the {\em non}-Fuchsian case we encounter very
different behaviours, both of solutions and of
mappings between  linear systems with such a singularity.
A generic solution of a non-Fuchsian system
\[y'=\frac{1}{x^m}B(x)y, \quad m\geq 2\]  does not have polynomial
growth in sectors, and generic formal power series solutions of such a system
 (as well as formal equivalences between generic non-Fuchsian
systems) are  divergent. The dynamics associated to a non-Fuchsian
singularity correspond to the dynamics of the
 vector field $$x^m\frac{\partial}{\partial x}+A(x)y\frac{\partial}{\partial y},$$
which is always resonant, in the sense of Poincar\'e-Dulac.

Further information on the classification of isolated
singularities  can be found
in e.g. \cite{MR2363178} or \cite{vazow}.

  Fuchsianity  admits a certain extention to the non-linear case as well,
giving rise to the notion of {\em Briot-Bouquet type ODEs}, that is, ODEs of
the form
\begin{equation}\Label{bbODE}\label{bbODE}  xy'=F(x,y), \end{equation}
where $x$ lies
in a neighborhood of $0$ in $\CC{}$, $y$ is $n$-dimensional and $F$ is
holomorphic in a neighborhood of $0$ in $\CC{n+1}$. Briot-Bouquet ODEs are
similar to linear systems of ODEs with a Fuchsian singularity
in many respects; for example, their formal power
series solutions are necessarily {\em convergent} (see, e.g., \cite{laine}).
Dynamics associated to a Briot-Bouquet type ODE corresponds to the dynamics of
the vector field $$x\frac{\partial}{\partial x}+F(x,y)\frac{\partial}{\partial
y}.$$

We also note that a Briot-Bouquet type ODE whose  principal matrix $F_y(0,0)$
has no positive integer eigenvalues  has at least one
holomorphic solution (see \cite{laine}).

\section{The associated ODE approach to the mapping problem}
We consider a real-analytic hypersurface with defining
equation as in  \eqref{madmissiblereal}.
The complex defining function of such a hypersurface
is given by
\begin{equation}\Label{madmissiblecomp}\label{madmissiblecomp}
w=\bar w+i \bar w^m\left(\epsilon |z|^2+\sum_{k,\ell\geq 2}\Theta_{k\ell}(\bar w)z^k\bar z^\ell\right).
\end{equation}
We recall from \autoref{ss:nonminimal} that
this means that  the Segre family
 $\mathcal S=\{Q_{(\xi,\eta)}\}$ of $M$
 is given by:
\begin{equation}\Label{mc}\label{mc}
w=\bar\eta e^{i\bar\eta^{m-1}\varphi(z,\bar\xi,\bar\eta)},\quad\mbox{where}
\quad\varphi(z,\bar\xi,\bar\eta)=\epsilon z\bar\xi+\sum_{k,\ell\geq 2}\varphi_{k\ell}(\bar\eta)z^k{\bar\xi}^\ell
\end{equation}

We will need the following fact proved in \cite{ekl}:
\begin{lemma}\Label{specialmap}\label{specialmap}(see \cite{ekl}).
Let  $H(z,w) =  \bigl(F(z,w),G(z,w)\bigr)$ be a formal transformation
vanishing at the origin, with invertible Jacobian $H'(0)$, which maps  a hypersurface defined by \eqref{madmissiblereal}
or equivalently \eqref{mc} into another such hypersurface. Then
$H$ satisfies
\begin{equation}\Label{specialg}\label{specialg}
\begin{gathered}
F_z(0,0)=\lambda, \quad  G_w(0,0)=\mu, \quad  G= O(w), \\ G_z=O(w^{m+1}),\quad \mu^{1-m}=|\lambda|^2, \quad \lambda\in\CC{}\setminus\{0\}, \, \mu\in\RR{}.
\end{gathered}
\end{equation}
In addition, we have
\begin{equation}\Label{taureal}\label{taureal}
G_{w^\ell}(0,0)\in\RR{}, \quad \text{ for } \ell \leq m.
\end{equation}
\end{lemma}

\autoref{specialmap}  implies in particular
that any transformation $H$ between hypersurfaces defined
by equations of the form \eqref{madmissiblereal} can be factored as
$$H=H_0\circ \psi,$$
for some dilation $\psi$ of the form \eqref{dilations}
and where   $H_0$ is
a transformation of the form:
$$z\mapsto z+f(z,w), \quad w\mapsto w+wg_0(w)+w^{m}g(z,w)$$
with
\begin{equation}\Label{normalmap}\label{normalmap}
 f_z(0,0)=0,\quad   g_0(0)=0,\quad g(z,w)=O(zw), \quad g_0^{(\ell)}(0)\in\RR{}, \quad \ell \leq m-1.
\end{equation}
(for $m=1$ the last condition is void).
In fact, one can also represent $H$ as $H=\psi\circ H_0$
(with a different $H_0$).
We therefore consider the classification problem only under transformations \eqref{normalmap}.

We now recall that  \cite{nonminimalODE,analytic} showed that
we can associate to a hypersurface in the
form \eqref{madmissiblereal} a second order
singular holomorphic ODE $\mathcal E(M)$ given by
\beq\Label{ODE}\label{ODE}
w''=w^m\Phi\left(z,w,\frac{w'}{w^m}\right),
\eeq
where $\Phi(z,w,\zeta)$ is holomorphic near the origin in $\CC{3}$, and satisfies $\Phi=O(\zeta^2)$.
This ODE is characterized by the condition that
any of the functions $w (z) = \Theta (z, \xi,\eta)$,
for $(\xi,\eta ) \in \bar U$, is a solution of the
 ODE \eqref{ODE}). We will decompose  $\Phi$ as
\begin{equation}\Label{expandPhi1}\label{expandPhi1}
\Phi(z,w,\zeta)=\sum_{j,k\geq 0,\ell\geq2} \Phi_{jk\ell}z^kw^j\zeta^\ell
\end{equation}
or
\begin{equation}\label{expandphi}
\Phi(z,w,\zeta)=\sum_{k\geq 0,l\geq 2}\Phi_{kl}(w)z^k\zeta^l.
\end{equation}

We now recall the approach used in \cite{KLS} and \cite{ekl}. Considering the transformation rule for second order ODEs and adapting it to ODEs \eqref{ODE} and maps \eqref{normalmap} expanded as  $\tilde f(z,w) = z + f(z,w) $, and $\tilde g (z,w) = w + w g_0 (w) + w^m g(z,w)  $, we get (see \cite{KLS},\cite{ekl}): 
\begin{multline}\Label{trule3}\label{trule3}
\Phi\left(z,w,\zeta\right)=\frac{1}{J}\Bigl[\bigl(1+ f_z+w^m f_w\cdot\zeta)^3 (1+g_0(w)+w^{m-1}g\bigr)^m\cdot\\
\cdot \Phi^*\Bigl(z+ f,w+wg_0(w)+ w^mg,\frac{g_z+\zeta(1+ wg_0'+g_0+mw^{m-1}g+w^mg_w)}{(1+g_0(w)+w^{m-1}g)^m(1+f_z+w^m\zeta f_w)}\Bigr)+\\
+ I_0(z,w)+I_1(z,w)\zeta+I_2(z,w)w^{m}\zeta^2+I_3(z,w)w^{2m}\zeta^3\Bigr],
\end{multline}
where $\zeta:=\frac{w'}{w^m}$ and
\beq\Label{2jet}\label{2jet} \begin{aligned}
J&=(1+f_z)(1+wg_0'+g_0+mw^{m-1}g+w^mg_w)-w^mf_wg_z,\\
I_0 &= g_zf_{zz}- (1+f_z)g_{zz},\\
I_1 &=\bigl(1+wg_0'+g_0+mw^{m-1}g+w^mg_w\bigr)f_{zz}-w^mf_wg_{zz}-\\
&-2(1+ f_z)
(mw^{m-1}g_z+w^mg_{zw})+2w^m g_z f_{zw},\\
I_2 &=w^m g_{z}f_{ww}-(1+ f_z)(wg_0''+2g_0'+m(m-1)w^{m-2}g+2mw^{m-1}g_w+w^mg_{ww})-\\
&-2f_w(mw^{m-1}g_z+w^mg_{zw}) +2(1+wg_0'+g_0+mw^{m-1}g+w^mg_w) f_{zw},\\
I_3 &=(1+wg_0'+g_0+mw^{m-1}g+w^mg_w)f_{ww}-\\
&-f_w(wg_0''+2g_0'+m(m-1)w^{m-2}g+2mw^{m-1}g_w+w^mg_{ww}).
\end{aligned}\eeq
Importantly, \eqref{trule3} is an identity in the {\em free variables $z,w,\zeta$}, where the latter triple runs a suitable open neighborhood of the origin in $\CC{3}$.

We recall then that, by collecting in \eqref{trule3} terms with $z^kw^j\zeta^l,\,l=0,1$, we obtain a system of PDEs of the kind:
\begin{equation} \Label{system1}\label{system1}
\begin{aligned}
f_{zz}=U(z,w,g_0,g_0',f,g,f_z,g_z,f_w,g_w,f_{zw},g_{zw}), \\
g_{zz}=V(z,w,g_0,g_0',f,g,f_z,g_z,f_w,g_w,f_{zw},g_{zw})
\end{aligned}
\end{equation}
for some  germs of holomorphic functions  $U,V$ at the origin.
Given a choice of (respectively holomorphic or formal) data
$$f(0,w)=f_0(w), \quad f_z(0,w)=f_1(w), \quad g(0,w)=0, \quad g_1(0,w)=g_1(w), $$
the Cauchy-Kowalevskaya theorem
guarantees
the existence of a unique (respectively holomorphic or formal) solution to \eqref{system1} with this  data.

The associated functions $\tilde f (z,w) = z + f(z,w)$, $\tilde g (z,w) = w + w g_0 (w) + g(z,w)$ transform $\mathcal{E}^*$ to the (up to the
initial data unique) $\mathcal{E}$. The initial conditions also
imply that $(\tilde f, \tilde g)$ is
of the form required in \eqref{normalmap}. To determine then the Cauchy data
\begin{equation}\Label{Cdata}\label{Cdata}
Y(w):=\bigl(f_0(w),f_1(w),g_0(w),g_1(w)\bigr),
\end{equation}
we collect in \eqref{trule3} terms with $z^kw^j\zeta^l,\,j=0,1,\,l=2,3$. This gives us a system of {\em singular} second order ODEs: 
\begin{equation}\Label{merom}\label{merom}
\begin{aligned}
w^{m+1}g_0''&=T_1(w,g_0,g_1,f_0,f_1,wg_0',w^{m}g_1',w^{m}f_0',w^{m}f_1'),\\
w^{2m}g_1''&=T_2(w,g_0,g_1,f_0,f_1,wg_0',w^{m}g_1',w^{m}f_0',w^{m}f_1'),\\
w^{2m}f_0''&=T_3(w,g_0,g_1,f_0,f_1,wg_0',w^{m}g_1',w^{m}f_0',w^{m}f_1'),\\
w^{2m}f_1''&=T_4(w,g_0,g_1,f_0,f_1,wg_0',w^{m}g_1',w^{m}f_0',w^{m}f_1')
\end{aligned}
\end{equation}
(we again refer to \cite{KLS},\cite{ekl} for details). 

Our Fuchsian type condition is obtained by requiring that, roughly speaking, the arising system of ODEs \eqref{merom} is Fuchsian (Briot-Bouquet). This is explained in the next section

\section{Fuchsian type ODEs and regularity of formal mappings}

\subsection{The normal form problem for Fuchsian type hypersurfaces}
First, we translate the Fuchsian type condition for hypersurfaces \eqref{madmissiblereal} described in the Introduction onto the language of associated ODEs. For  the functions $\Phi,\Phi^*$, we make use of the expansion
\eqref{expandphi}. We now introduce
\begin{definition}\Label{fuchsian}\label{fuchsian}
 An ODE $\mathcal{E}$, defined by
  \eqref{ODE},
  is called {\em Fuchsian} (or {\em a Fuchsian type ODE}), if $\Phi$ satisfies the conditions:
\begin{equation}\Label{FODE}\label{FODE}
\begin{aligned}
&\ord\,\Phi_{02}(w)\geq m-1; \,\ord\,\Phi_{03}(w)\geq 2m-2; \,\ord\,\Phi_{12}(w)\geq m-1;\,\ord\,\Phi_{13}(w)\geq 2m-2;\\
&\ord\,\Phi_{0l}(w)\geq 2m-l+2, \,\,4\leq l\leq 2m+1;\,\,
\ord\,\Phi_{k2}(w)\geq 2m-k, \,\,2\leq k\leq 2m+1; \\
&\ord\,\Phi_{kl}(w)\geq 2m-k-l+3,\,\,k\geq 1,\,l\geq 3,\,5\leq k+l\leq 2m+2.
\end{aligned}
\end{equation}
\end{definition}
We make use of the following:
\begin{proposition}[See \cite{ekl}]\Label{transferfuchs}\label{transferfuchs}
For a Fuchsian type hypersurface $M\subset\CC{2}$, its associated ODE $\mathcal E(M)$ is of Fuchsian type as well.

\end{proposition}
We next prove the invariance of the Fuchsian type condition.
\begin{theorem}\Label{fuchsinv}\label{fuchsinv}
The property of being Fuchsian for a hypersurface \eqref{madmissiblereal} does not depend on the choice of (formal or holomophic) coordinates of the kind \eqref{madmissiblereal}.
\end{theorem}
\begin{proof}
In view of \autoref{fuchsian}, we can switch to associated ODEs and it is enough to prove the invariance of the Fuchsianity for them. As discussed above, we can restrict to transformations \eqref{normalmap}. Let us consider then the transformation rule \eqref{trule3} (with a fixed transformation within it), when the source ODE (with the defining function $\Phi^*$) is of Fuchsian type. We then claim the following: for all the coefficient functions $\Phi_{kl},\,\,k\geq 0,\,\,l\geq 2$ involved in the Fuchsianity conditions \eqref{FODE}, {\em with the exception of the coefficients functions} $\Phi_{k2},\Phi_{k2}^*,\,\,k\geq 2$, the Fuchsian conditions \eqref{FODE} are satisfied.
 Indeed, we fix any $(k,l)$ relevant to \eqref{FODE}, and from the transformation rule \eqref{trule3} we can see that the target coefficient function $\Phi_{kl}$ is a sum of three groups of terms: (i) terms $\Phi^*_{\alpha\beta}$ with $\alpha+\beta\geq k+l$ which are multiplied by a power series in $w$ with order at $0$ at least $k+l-\alpha-\beta$; (ii) terms $\Phi^*_{\alpha\beta}$ with $\alpha+\beta< k+l$; (iii) terms arising from the expressions $I_j,\,0\leq j\leq 3$ (relevant for $l=2,3$ only). In view of the linearity of the Fuchsianity conditions in $k,l$, it is not difficult to see that terms of the first kind all have order at $0$ at least as the one required for the Fuchsianity. Terms of the second kind already all have order bigger than the one required for Fuchsianity. Finally, terms of the third kind automatically provide order at least $2m$ sufficient for the Fuchsianity, except for the case $l=2$. For $k=0,1$ and $l=2$ though even the automatically provided order $m$ suffies, and this proves the claim.

It remains to deal with terms $\Phi_{k2}$ with $k\geq 2,\,2\leq k\leq 2m+1$. We note, however, that the ODEs under consideration have a real structure, which is why (in view of the reality condition) we have
\begin{equation}\Label{heq}\label{heq}
\ord\,h_{kl}(w)=\ord\,h_{lk}(w)
\end{equation}
for all $k,l$. This, in view of the transfer relations between $\Phi$ and $h$, gives, in particular:
$$\ord\,\Phi_{k2}(w)=\ord\,h_{k+2,2}(w)=\ord\,h_{2,k+2}(w)=\ord\,\Phi_{0,k+2}(w)\geq 2m-k$$
(the last inequality follows from the Fuchsianity condition for $\Phi_{0,k+2}$ being already proved). This finally proves the theorem.

\end{proof}

We now proceed with the proof of \autoref{theor5}. We follow the scheme in Section 3, and obtain a system of singular ODEs of the kind \eqref{merom} for the Cauchy data $Y(w)$, as in\eqref{Cdata}, assuming the source ODE (with the defining function $\Phi^*$) is of Fuchsian type. For the purposes of this section, we prefer to write down the obtained system in the form
\begin{equation}\Label{meromF}\label{meromF}
w^{m+1}g_0''=S\bigl(w,Y(w),wY'(w)\bigr),\,\, w^{2m}X''=T\bigl(w,Y(w),wY'(w)\bigr),
\end{equation}
where $$X(w):=(g_1(w),f_0(w),f_1(w)),\quad Y(w):=(g_0(w),X(w)),$$
and $S,T$ are holomorphic near the origin.

For the functions $T,S$ we will use the expansion
\begin{equation}\Label{expandTS}\label{expandTS}
T(w,Y,\tilde Y)=\sum_{\alpha,\beta\geq 0}T_{\alpha,\beta}(w)Y^\alpha\tilde Y^\beta,
\end{equation}
where $\alpha,\beta$ are multiindices, and similarly for $S$.   We now shall prove the following key
\begin{proposition}\Label{Tkl}\label{Tkl}
Under the Fuchsian type condition, the coefficient functions $T_{\alpha,\beta}(w),S_{\alpha,\beta}(w)$ satisfy
\begin{equation}\Label{Tkle}\label{Tkle}
\ord T_{\alpha,\beta}\geq 2m-1-|\alpha|-|\beta|,\,\,\ord S_{\alpha,\beta}\geq m-|\alpha|-|\beta|,\,\,|\alpha|+|\beta|>0.
\end{equation}
\end{proposition}
\begin{proof}
For the proof, we make use of \eqref{FODE} (applied for the source defining function $\Phi^*$), and then study carefully  the contribution of terms $\Phi^*_{kl}$ into the basic identity \eqref{trule3}. Let us fix for the moment some positive value of $|\alpha|+|\beta|$. Then it is straightforward to check, by considering  \eqref{trule3}, that $T_{\alpha,\beta}$ as above can arise only from $\Phi^*_{kl}$ with $k+l\leq  |\alpha|+|\beta|+4,$ while $S_{\alpha,\beta}$ as above can arise only from $\Phi^*_{kl}$ with $k+l\leq  |\alpha|+|\beta|+2$. (And in the latter cases a respective $\Phi_{kl}^*$ is a factor for $Y^\alpha(wY')^\beta$). Now it is not difficult to verify that \eqref{FODE} implies \eqref{Tkl}.
\end{proof}

\begin{corollary}\Label{lowterms}\label{lowterms}
For the $(0,0)$ coefficient functions in \eqref{meromF} we have
\begin{equation}\Label{00}\label{00}
\ord\,S_{0,0}\geq m;\quad \ord\,T_{0,0}\geq 2m-1.
\end{equation}
As a consequence, for the target ODE defining function $\Phi$ we have:
\begin{equation}\Label{Phi23a}\label{Phi23a}
\begin{aligned}
\Phi_{0j2}=0,&\,\,0\leq j\leq m-2; \\
 \Phi_{1j2}=\Phi_{0j3}=\Phi_{1j3}=0,&\,\,0\leq j\leq 2m-3;\\
\Phi_{0,m-1,2}=\Phi^*_{0,m-1,2};\,\,\Phi_{0,2m-2,3}=\Phi^*_{0,2m-2,3};&\,\,\Phi_{1,2m-2,2}=\Phi^*_{1,2m-2,2};\,\,\Phi_{1,2m-2,3}=\Phi^*_{1,2m-2,3}.
\end{aligned}
\end{equation}
\end{corollary}
\begin{proof}
As follows from the definition of $S_{\alpha,\beta},T_{\alpha,\beta}$ and the Fuchsianity, all terms in the first equation in \eqref{meromF} have order at least $m$ in $w$ with possibly the exception of terms arising from $S_{0,0}$, while all terms in the second equation in \eqref{meromF} have order at least $2m-1$ in $w$ with possibly the exception of terms arising from $T_{0,0}$. This proves \eqref{00}. To prove \eqref{Phi23a}, we note that the $(m-1)$-jet of $S_{0,0}$ and the $(2m-2)$-jet of $T_{0,0}$ respectively are formed from differences  between coefficients $\Phi_{kjl}$ and $\Phi^*_{kjl}$ aparent in \eqref{Phi23a}, and this proves \eqref{Phi23a}.
\end{proof}

We shall now prove that any solution of the system of singular ODEs \eqref{meromF}. In view of the discussion in Section 3, this would imply the convergence of the formal map between the given ODEs \eqref{ODE} and  the given real hypersurfaces, and hence the assertion of \autoref{theor5}. 

Let $H(w)$ be such a formal solution of \eqref{meromF}. We decompose it as
\begin{equation}\Label{trick}\label{trick}
H(w)=P(w)+Z(w),
\end{equation}
where $P(w)$ is a  polynomial without constant term of degree $\leq 2m-1$, 
while where $Z(w)$ is a formal series
of the kind $O(w^{2m})$.
The substitution \eqref{trick} (for a fixed $(P(w)$) turns \eqref{meromF} into a similar system of ODEs for the unknown function $Z(w)$. We shall now prove
\begin{lemma}\Label{0110}\label{0110}
The transformed  system (in the same way as the initial system)  satisfies
\begin{equation}\Label{0110a}\label{0110a}
\ord\,\tilde S_{01}\geq m-1,\,\,\ord\,\tilde S_{10}\geq m-1,\,\,\ord\,\tilde T_{01}\geq 2m-2,\,\, \ord\,\tilde T_{10}\geq 2m-2
\end{equation}
(the tilde here stands for coefficients of the transformed system).
\end{lemma}
\begin{proof}
The proof of the lemma is obtained by putting together the expansion \eqref{expandTS}, the conditions \eqref{Tkl}, and the fact that $P(w)$ is vanishing at the origin.
\end{proof}
Now, based on \autoref{0110}, we perform the substitution
\begin{equation}\Label{ZtoU}\label{ZtoU}
Z:=w^{2m}U,
\end{equation}
which turns the "tilde" system into a new system of four meromorphic ODEs for the unknown function $U$, which, according to \eqref{trick}, has a formal solution $U(w)$ vanishing at the origin. It is straightforward to check then, by combining \eqref{ZtoU} and \eqref{0110a}, that the new system system can be written in the form
\begin{equation}\Label{bb2}\label{bb2}
w^2U'=R(w,U,wU'),
\end{equation}
where $R$ is a holomorphic  function defined near the origin. Performing finally in the standard fashion the substitution
$$V:=wU'$$
and introducing the extended vector function $\bold U:=(U,V)$, we obtain a first order ODE
\begin{equation}\Label{bb1}\label{bb1}
w\bold U'=Q(w,\bold U'),
\end{equation}
where $Q$ is a holomorphic near the origin function. The ODE \eqref{bb1} is a Briot-Bouquet type ODE (see Section 2), hence its formal solutions are convergent, as required.

This completes the proof of \autoref{theor5}. \qed

\section{Regularity of smooth mappings between Fuchsian type hypersurfaces} In this section we shall prove \autoref{theor4}.
Compared to the proof of \autoref{theor5}, we need an additional argument, which is the following regularity result for  Fuchsian (Briot-Bouquet)  systems of meromorphic ODEs.
\begin{proposition}\Label{analyticity}\label{analyticity}
Consider a first order real ODE
\begin{equation}\Label{1order}\label{1order}
xy'=F(x,y),\quad x\in[0,a],
\end{equation}
with $y$ being $n$-dimensional, $n\geq 1$, and $F$ analytic.  Assume it has a solution $y(x)$ which is $C^\infty$ on $[0,a]$. Then $y(x)$ is analytic everywhere on $[0,a]$.
\end{proposition}

\begin{remark}\Label{bb}\label{bb}
A singular ODE  \eqref{1order} belongs to the classical class of {Briot-Bouquet} type ODEs discussed in Section 2. Their formal solutions at the singular point $x=0$ are {\em convergent}, which, however does {\em not} say anything about the regularity of smooth solutions, which is why \autoref{analyticity} requires a separate proof.
\end{remark}
\begin{proof}[Proof of \autoref{analyticity}]
The analyticity of $y(x)$ everywhere outside $x=0$ follows from the analyticity of the given ODE, which is why we consider only the analyticity at the singularity $x=0$. First, consider the Taylor series $\hat y(x)$ of $y(x)$. Since, again, \eqref{1order} is a Briot-Bouquet ODE, $\hat y(x)$ is convergent. Hence, taking $y-\hat y(x)$ as a new unknown function, we get an ODE again of the kind \eqref{1order}  which has now a {\em flat} at $x=0$ solution on $[0,a]$. We assume, by contradiction, that this solution is not identical zero near $x=0$. Substituting the flat solution into the new ODE\eqref{1order} and equalizing the Taylor series in both sides, we conclude that $F(x,0)=0$. Hence we conclude that the (again analytic) right hand side expands as
$$F(x,y)=A(x)y+\cdots,$$
where $A(x)$ is an analytic at the origin matrix, and dots stand for terms of degree at least $2$ in $y$.

Second, let us use the notation $|y(t)|$ for the Euclidean norm, and $||y||$ for the sup norm of $y$ on $[0,a]$. Since $y$ is flat at $0$, we may shrink the interval to make $\|y\|$ small. Using the analyticity of $F$, we then have the bound
\begin{equation}\Label{bound1}\label{bound1}
|F(x,y(x))|\leq C|y(x)|,
\end{equation}
where $C$ is a constant depending on $\|y\|$.

Third, we make a simple observation that $|y|$ can not vanish for $x>0$. Indeed, any solution with $y(x_0)=0,\,x_0\neq 0$ would need to be identical zero by uniqueness near $x_0$, and hence identical zero by the analyticity of the ODE.

Fourth, we do the following: we "resolve the singularity" of \eqref{1order} by making the substitution
$$x:=e^t,\quad t\in(-\infty,\ln a].$$
Now the ODE \eqref{1order} reads as
\begin{equation}\Label{1ord}\label{1ord}
\frac{dy}{dt}=F(e^t,y)=:\tilde F(t,y).
\end{equation}
We denote the new solution by $y(t)$ and still have
\begin{equation}\Label{bound2}\label{bound2}
|\tilde F(t,y(t))|\leq C|y(t)|.
\end{equation}

Now we need to obtain certain bounds. Taking the limit in the triangle inequality, we have
$$\frac{d}{dt}|y(t)|\leq \left|\frac{dy}{dt}\right|.$$
In view of this and the inequality \eqref{bound2},
$$\frac{d}{dt}\ln|y(t)|=\frac{1}{|y(t)|}\frac{d}{dt}|y(t)|\leq \frac{1}{|y(t)|}\left|\frac{dy}{dt}\right|\leq C,$$
and by integrating over $[t,\ln a]$ we obtain:
\begin{equation}\Label{bound3}\label{bound3}
\ln |y(\ln a)|-\ln|y(t)|\leq C(\ln a-t).
\end{equation}
Simplifying \eqref{bound3} and applying exp, we finally get for the initial function $y(x)$:
\begin{equation}\Label{bound4}\label{bound4}
|y(x)|\geq \tilde C\cdot x^C
\end{equation}
($\tilde C$ is some other constant, which is nonzero since $|y(a)|$ is nonzero!). But \eqref{bound4} is a contradiction with the fact that $y(x)$ is flat near $0$, and this proves the desired analyticity statement.
\end{proof}
\begin{remark}\Label{cxbb}\label{cxbb}
The assertion of \autoref{analyticity} holds also for a {\em complex} Briot-Bouquet ODE, i.e. when $y(x)$ is complex-valued and $F$ is complex analytic (one just has to split the real and imaginary parts, and this immediately gives an already {\em real} ODE \eqref{1order} for the vector function formed from the real and imaginary parts of $y$).
\end{remark}
We are now in the position to prove \autoref{theor4}.
\begin{proof}[Proof of \autoref{theor4}]
We come back to the proof of \autoref{theor4}. Note that a hypersurface \eqref{madmissiblereal} necessarily contains  (the germ at the origin of) the real line $L=\{z=0,\,\,\im w=0\}$. This means, in particular, that for the given map $H(z,w)$,  the vector functions $H(0,w),H_z(0,w)$ are well defined on $L$ and are holomorphic in its open neighborhood. Arguing now identically to the above proof of \autoref{theor5}, we reduce the analyticity problem for the given CR-map to the analyticity of $C^\infty$ smooth solutions of an ODE identical to \eqref{bb1}. The only difference is that, instead of substituting a formal power series map into the basic identity \eqref{trule3}, we substitute into \eqref{trule3} a holomorphic map in a  domain $\Omega$, containing $0$ in its closure and coming from the analyticity of the map in a neighborhood of the Levi-nondegenerate part of $M$. In view of the above, the Cauchy data \eqref{Cdata} of the map $H$ is $C^\infty$ on the real line, and so is a solution of \eqref{bb1} under discussion. We then apply \autoref{analyticity} (together with \autoref{cxbb}) and conclude that the desired solution of \eqref{bb1} is analytic, and so is the Cauchy data \eqref{Cdata} and hence the map $H$. This completely proves the theorem.

\end{proof}

\bibliographystyle{plain}
 \bibliography{bibfile}









\end{document}